\documentclass[11pt]{article}
 \usepackage{graphicx,amsmath,amsfonts,amssymb}
\usepackage{epsfig}

\def\mbf#1{\mbox{\boldmath $#1$}}

 \newcommand{\dd}[1]{\if#1N I\hspace{-1.2mm}N \else {\if#1R I\!\!R \else
{\if #11 1\!\!1 \else #1\!\!\!#1 \fi} \fi} \fi}

\newcommand{\eqnb}{\begin{equation}}
\newcommand{\eqne}{\end{equation}}
\newcommand{\Proof}{{\bf Proof: }}

\newcommand{\p}{{\mathbb P}}
\newcommand{\e}{{\mathbb E}}

\newcommand{\be}{\mbox{\boldmath $e$}}
\newcommand{\bpi}{\mbox{\boldmath $\pi$}}

\newtheorem{Thm}{Theorem}[section]

\newtheorem{Def}{Definition}[section]

\newtheorem{Cor}{Corollary}[section]
\newtheorem{Lem}{Lemma}[section]

\newtheorem{remark}{Remark}[section]
\newtheorem{Cond}{Condition}[section]

\def\S{$\hfill\square$}

\title{\Large \bf Ergodicity for the $GI/G/1$-type Markov Chain\thanks{This research was supported in part by NSERC
Discovery Grants and an NSFC grant (10121101,10301007). The authors thank Dr.~Don Dawson for his comments and suggestions made throughout the research process
leading to this manuscript. The second author would also like to express his appreciations to Dr.~Dawson for his continuous encouragement and guidance during his studies
for his Ph.D. degree.}}
\author{YongHua Mao\\
School of Mathematical Science, Beijing Normal University\\
Beijing, 100875, P.R. China\\
\and Yongming Tai \\
School of Mathematics and Statistics, Carleton University \\
Ottawa, ON Canada K1S 5B6
\and Yiqiang Q. Zhao \\
School of Mathematics and Statistics, Carleton University \\
Ottawa, ON Canada K1S 5B6 \and Jiezhong Zou \\
School of Mathematical Science and Computing Technology \\
Central South University \\
Changsha, Hunan, 410075, P.R. China }
\date{\today}

\begin{document}
 \maketitle

 \begin{abstract}
Ergodicity is a fundamental issue for a stochastic process. In
this paper, we refine results on ergodicity for a general type of
Markov chain to a specific type  or the $GI/G/1$-type Markov
chain, which has many interesting and important applications in
various areas. It is of interest to obtain conditions in terms of
system parameters or the given information about the process,
under which the chain has various ergodic properties.
Specifically, we provide necessary and sufficient conditions for
geometric, strong and polynomial ergodicity, respectively.
 \end{abstract}

 \begin{minipage}{0.9\linewidth}
  \end{minipage}

\section{Introduction}

Ergodicity is a fundamental issue in the study of a stochastic process. There are many references in this area, among which closely
related to our study include books by Anderson~\cite{Andersonmc}, Meyn and Tweedie~\cite{MeynStability}, Chen~\cite{chen:05}, and references
therein. Results in these books are usually presented for a general type of stochastic process, or Markov chain, while our focus is on a specific type of Markov
chain; that is, the $GI/G/1$-type Markov chain. We apply standard methods used for a general process, combined with various techniques in dealing with
block-structured matrices, to characterize ergodic properties in terms of system parameters, which is an extension of the existing research.

The $GI/G/1$-type Markov chain is a very important type of
block-structured stochastic process with many applications in
queueing theory; for example, see Grassmann and
Heyman~\cite{Grassmann:14}, Zhao, Li and
Braun~\cite{Zhao:11,fact-mcap}, and Zhao~\cite{Zhao:26}. We refine
literature results on general Markov chains to obtain ergodicity
conditions for this specific type of Markov chain. In this paper,
we focus on the study of ergodicity of this type of Markov chain.
For the $GI/G/1$-type Markov chain, the ordinary ergodicity has
been well studied. Necessary and sufficient conditions have been
reported in the literature; for example, see
Asmussen~\cite{Asmussen:07}, Zhao, Li and
Braun~\cite{Zhao:11,fact-mcap}, and Zhao~\cite{Zhao:26}. This will
not be discussed again in the current paper. Instead, we will
consider three other types of ergodicity: geometric, strong (or
uniform) and polynomial. Related results were reported in Spieksma
and Tweedie~\cite{st}.

Related to our research, H{\o}jgaard and
M{\o}ller~\cite{Hojgaard1996}, using the coupling method and
stopped random walks, provided a sufficient condition for
geometric and polynomial ergodicity, respectively, for the
$GI/G/1$-type Markov chain. Hou and Liu~\cite{HouEergodicity}
derived a necessary and sufficient condition for polynomial
ergodicity for the $M/G/1$ queue by analyzing the generating
function of the first return probability, and extended their study
to the $M/G/1$-type Markov chain in Liu and Hou~\cite{Liu:35}.
Jarner and Tweedie~\cite{Jarner2001} proved that for
random-walk-type Markov chains, the geometric (light) and
polynomial tail asymptotics in the stationary probability
distribution are necessary  for the geometric ergodicity and
polynomial ergodicity, respectively.

The main contributions in this paper include necessary and
sufficient conditions (Theorem~\ref{thm:equiv},
Theorem~\ref{thm:strong}, Theorem~\ref{thmergodic02} and
Theorem~\ref{thmergodi02}) for each of these three types of
ergodicity, given in terms of system parameters, or the given
information about the $GI/G/1$-type Markov chain.

Ergodicity and tail asymptotics of the stationary probability
distribution are usually studied separately because of the obvious
distinction between these two concepts. For a stochastic process,
ergodicity deals with conditions under which the marginal
distribution at time $t$ converges to its limiting distribution in
various speeds as $t$ goes to infinity, while the tail asymptotic
of the stationary (limiting) probability distribution is concerned
with the speed to zero of the tail probabilities. It is
interesting to observe that the same necessary and sufficient
condition for both geometric ergodicity and a geometric decay in
the stationary probability distribution immediately allows us to
draw the conclusion of ergodicity on a number of important models
for which the tail asymptotics have been known, and vice versa.
The equivalence also opens a new door for us to take advantage of
possibly using newly developed approaches and known results in
studying ergodicity to study tail asymptotics of the stationary
probability distribution. For a symmetric Metropilis-Hastings
algorithm, Mengersen and Tweedie \cite{met} and Jarner and Hansen
\cite{jh} proved that these two concepts are actually equivalent.
It is our goal to provide necessary and sufficient conditions for
ergodicity for the $GI/G/1$-type model in this paper.

The rest of the paper is organized into four sections. In
Section~2, the $GI/G/1$-type Markov chain is reviewed, and a
spectral property, Lemma \ref{lem:eigen}, is obtained that plays a
key role in proving the main result for geometric ergodicity.
Sections~3--5 study geometric, strong and polynomial ergodicity,
respectively.

\section{The $GI/G/1$-type Markov chain}

Consider a discrete time irreducible aperiodic Markov chain, whose
transition probability matrix is given by
 \begin{eqnarray} \label{eqn:P}
  P = \left(\begin{array}{ccccc}
   P_{0,0} & P_{0,1} & P_{0,2} & P_{0,3} & \cdots\\
   P_{1,0} & P_{1,1} & P_{1,2} & P_{1,3} & \cdots\\
   P_{2,0} & P_{2,1} & P_{2,2} & P_{2,3} & \cdots\\
   P_{3,0} & P_{3,1} & P_{3,2} & P_{3,3} & \cdots\\
   \vdots & \vdots & \vdots & \vdots & \cdots
  \end{array}\right)
  = \left(\begin{array}{ccccc}
   B_0 & B_1 & B_2 & B_3 & \cdots\\
   B_{-1} & A_{0} & A_{1} & A_{2} & \cdots\\
   B_{-2} & A_{-1} & A_{0} & A_{1} & \cdots\\
   B_{-3} & A_{-2} & A_{-1} & A_{0} & \cdots\\
   \vdots & \vdots & \vdots & \vdots & \cdots
  \end{array}\right),
 \end{eqnarray}
where $A_{i}$ and $B_{i}$ for $i = 0, \pm 1, \pm 2, \cdots$, are
matrices of size $m\times m$. The state space for the Markov chain
$P$ can be expressed by $S=\cup_{i=0}^{\infty}L_{i}$, where
$L_{i}=\{(i,j);j=1,2,\ldots,m \}$ for $i \geq 0$. In a state
$(i,j)$, $i$ is referred to as a {\em level} and $j$ as a {\em
phase}. We also write $L_{\leq i} =\cup_{k=0}^{i}L_{k}$.
\begin{remark}
In fact, we can more generally assume that $B_0$ is a matrix of
size $m_0 \times m_0$ with $m_{0}\neq m$ for the results obtained
in this paper.
\end{remark}

Along the same line as in~\cite{Zhao:11}, we define the
$R$-measures $R_{i,j}$ for $i<j$ and the $G$-measures $G_{i,j}$
for $i>j$ for the $GI/G/1$-type Markov chain. $R_{i,j}$ is a
matrix of size $m\times m$ whose $(r,s)$th entry is the expected
number of visits to state $(j,s)$ before hitting any state in
$L_{\leq (j-1)}$, given that the process starts in state $(i,r)$.
$G_{i,j}$ is a matrix of size $m\times m$ whose $(r,s)$th entry is
the probability of hitting state $(j,s)$ when the process enters
$L_{\leq (i-1)}$ for the first time, given that the process starts
in state $(i,r)$. We refer to the matrices $R_{i,j}$ and $G_{i,j}$
as the matrices of the expected number of visits to higher levels
before returning to lower levels and the matrices of the first
passage probabilities to lower levels, respectively.
From~\cite{Zhao:11}, we can write $R_{n-i} = R_{i,n}$ and $G_{n-i}
= G_{n,i}$ for $i>0$ due to the property of repeating rows. If the
$GI/G/1$-type Markov chain is positive recurrent, then the
stationary distribution $\{\bpi_{k}\}$ can be expressed in terms
of the $R$-measures (see~\cite{Grassmann:14}):
\eqnb\label{SE:0102}
   \bpi_{n} = \bpi_{0}R_{0,n} + \sum\limits_{k=1}^{n-1} \bpi_{k}R_{n-k}, \;\;\; n\geq 1.
\eqne
Define the generating functions for the stationary
distribution $\{\bpi_{k}\}$, the matrix sequences $\{R_{0,k}\} $
and $\{R_{k}\} $, respectively, as $\bpi^{*}(z) =
\sum_{k=0}^{\infty} z^{k}\bpi_{k}$, $R^{*}_{0}(z) =
\sum_{k=1}^{\infty} z^{k}R_{0, k}$ and $R^{*}(z) =
\sum_{k=1}^{\infty} z^{k}R_{k}$. Then, we have that
\eqnb\label{SE:0103}
   \bpi^{*}(z)[I - R^{*}(z)] = \bpi_{0}R^{*}_{0}(z).
\eqne

Throughout this paper, we assume that the phase process $A =
\sum_{k=-\infty}^{\infty} A_k$ is irreducible. Therefore, when $A$
is stochastic in addition, there is a unique invariant probability
vector $\mbf\mu$ of $A$, that is $\mbf\mu=\mbf\mu A$.

\begin{Def}
\label{Def:light} For a sequence $\left\{  c_{k}\right\}  $ of
nonnegative scalars, it is called light-tailed if
\[
\sum\limits_{k=1}^{\infty}c_{k} e^{\epsilon k} <+ \infty
\]
for some $\varepsilon>0$. For a sequence $\left\{  C_{k}\right\} $
of nonnegative matrices of size $m\times n$, it is called
light-tailed if for all $i=1,2,\cdots m$ and $j=1,2,\cdots,n$, the
sequences $\left\{ C_{k}\left(  i,j\right)  \right\}  $ of
nonnegative scalars are light-tailed, where $C_{k}\left(
i,j\right)  $ is the $\left( i,j\right)$th entry of $C_{k}$.
\end{Def}

Let $A^*_+(z;r,s)$ and  $B^*_+(z;r,s)$ be the $(r,s)$th entry of the matrix generating
function $A^*_+(z)=\sum_{k=1}^{\infty} A_k z^k$  and
$B^*_+(z)=\sum_{k=1}^{\infty} B_k z^k$, respectively. Denote by $\phi_{A+}(r,s)$ and $\phi_{B+}(r,s)$ the
convergence radii of $A(z;r,s)$ and  $B(z;r,s)$, respectively. Let
\[
    \phi_{A+} = \min_{r,s} \phi_{A+}(r,s)
\]
and
\[
    \phi_{B+} = \min_{r,s} \phi_{B+}(r,s).
\]
In Li and Zhao \cite{liz1, liz2}, they proved the following.

\begin{Lem}
\label{lem:light} When $A$ is stochastic, for the $GI/G/1$-type Markov chain, the
sequence of the stationary probability vectors $\{\mbf\pi_k\}_{k\ge1}$ is
light-tailed if and only if $\min(\phi_{A+}, \phi_{B+}) >1$, or
both $\{A_k\}_{k \geq 1}$ and $\{B_k\}_{k \geq 1}$ are
light-tailed.
\end{Lem}

Define $A^*(z) = \sum_{k=-\infty}^{\infty} z^k A_k$. Then,
according to Li and Zhao~\cite{liz2}, we have the following
spectral property.
\begin{Lem}
Suppose that the $GI/G/1$-type Markov chain defined in (\ref{eqn:P}) is irreducible, aperiodic and positive recurrent.
If $\Omega=\{1<|z|<\phi_{A+}: {\rm det}(I-A^*(z))=0\}$
is not empty, then there must exit a positive $z_0\in\Omega$ such that $z_0\le|z| $ for $z\in\Omega$.
\end{Lem}

Define $\eta=z_0$ for $\Omega\not=\emptyset$; $\eta=\infty$ for $\Omega=\emptyset$.
Let $\chi(z)$ for $z>0$ be the largest eigenvalue of $A^*(z)$.
 We have the following spectral property.

\begin{Lem}
\label{lem:eigen} Assume that $A$ is stochastic and $\phi_{A+} >1$ for the
irreducible aperiodic positive recurrent $GI/G/1$-type Markov chain defined in (\ref{eqn:P}). Then
$\chi(z)<1$  for any $1<z < \min(\phi_{A+},\eta)$,
and there exists a Perron-Frobenius eigenvector $\mbf Y(z)$ with $\mbf Y(z)\ge\mbf{e}^{t}$
such that
\begin{equation}\label{z1}
A^*(z)\mbf Y(z)=\chi(z)\mbf Y(z),
\end{equation}
where $\mbf e$ is a row vector of ones and $\mbf {e}^{t}$ is the transpose of $\mbf e$.
\end{Lem}

\Proof
Without loss of generality, assume that $\eta<\phi_{A+}$. Since $\chi(1)=\chi(\eta)=1$,
it follows from the continuity of $\chi(\cdot)$ and the definition of $\eta$
that we only need to show that $\chi'(1)<0$.

Differentiating both sides of equation (\ref{z1}) at $z=1$ gives
\[
{A^*}'(1)\chi(1)+{A^*}(1)\mbf Y'(1)=\chi'(1)\mbf Y(1)+\chi(1)\mbf Y'(1).
\]
Note that because $A^*(1)=A,
{A^*}'(1)=\sum_{k=-\infty}^{\infty}kA_k$ and $\mbf Y(1)=\mbf
{e}^{t}$, the lemma follows from the ergodicity condition
\[
\chi'(1)=\mbf\mu\sum_{k=-\infty}^{\infty}kA_k\mbf {e}^{t}<0.
\]
\S

\section{Geometric ergodicity}

In this section, we present a necessary and sufficient condition
for geometric ergodicity, which is the same condition for a
geometric tail in the stationary probability distribution.

Let $P^n((i,r),(j,s))$ be the $n$-step transition probability for
the irreducible and aperiodic Markov chain of $GI/G/1$ type given
in (\ref{eqn:P}). $P$ is called geometrically ergodic if
 there exists a rate $\rho<1$  such that
\[
|P^n((i,r),(j,s))- \pi_{j,s}|\le M_{(i,r)(j,s)}\rho^n, \quad \text{for all $n\ge0$ and for all $i,j,r,s$},
\]
where $\pi_{j,s}$ is the stationary probability in level $j$ and phase $s$ and $M_{(i,r)(j,s)} < \infty$. It follows from Theorem 4.31 in Chen \cite{cbk}
that geometric ergodicity is equivalent to the following
condition.

\begin{Cond} \label{con:1}
There exist a finite set $H\not=\emptyset$, a constant $\lambda<1$
and a finite function $V \geq 1$ defined on the state space of the Markov
chain $P$, such that
\begin{equation}\label{c1}
  \begin{cases}
    P V(i,r) \leq \lambda V(i,r), & \text{\rm for $(i,r)\not\in H$}, \\
    P V(i,r) <\infty, & \text{\rm for $(i,r)\in H$}.
  \end{cases}
\end{equation}

\end{Cond}
See also Anderson \cite{Andersonmc} or Meyn and Tweedie \cite{MeynStability}.

When $A$ is stochastic, the $GI/G/1$-type Markov chain is a
special case of the random-walk-type Markov chain in Jarner and
Tweedie~\cite{Jarner2001}, so it follows from Theorem 2.2 of
\cite{Jarner2001} that the condition in Lemma \ref{lem:light} or
 $\min(\phi_{A+}, \phi_{B+}) >1$
 is necessary for $P$ to be geometrically
ergodic.  To show the equivalence between the geometric ergodicity
and geometric stationary tail as stated in Lemma~\ref{lem:light},
we only need to show that under the assumption $\min(\phi_{A+},
\phi_{B+}) >1$, (\ref{c1}) holds.

\begin{Thm}\label{thm:equiv}
Assume that $A$ is stochastic and $\min (\phi_{A+}, \phi_{B+}) >1$ for the irreducible
aperiodic positive recurrent $GI/G/1$-type Markov chain defined
in (\ref{eqn:P}). (\ref{c1}) holds for $V(i,r) =
 z^{i} y_r$ for any $1<z <\min(\phi_{A+},\phi_{B+}, \eta)$,
 where $\mbf{Y} :=\mbf Y(z)= (y_1, y_2, \cdots, y_m)^t$ is the
 right eigenvector of   $A^*(z)$ corresponding to the eigenvalue $\chi(z)$
 given in Lemma~\ref{lem:eigen}.
\end{Thm}

\Proof
For any fixed $1<z <\min(\phi_{A+},\phi_{B+}, \eta)$, $\delta:=1-\chi(z)>0$ by Lemma \ref{lem:eigen}.
Let
$N$ be large enough such that $\alpha z^{-N}\le \delta/2$ or $N\ge (\log 2\alpha-\log\delta)/\log z$,
where $\alpha=\max_{1\le r\le m}y_r\ge1$.

Set $H=L_{\le N}$, a finite set. Then, for $(i,r)\not\in H$ or $i> N$, we have
\begin{equation} \label{eqn:con}
\aligned
   P\mbf v(i)&:=\sum_{k=0}^\infty P_{i,k}\mbf v_k
   = B_{-i} \mbf{v_0}  + A_{-i+1} \mbf{v_1} + A_{-i+2} \mbf{v_2} +
    A_{-i+3} \mbf{v_3}  +    \cdots  \\
    &=\left(B_{-i} z^{-i} + A_{-i+1} z^{-i+1}  + A_{-i+2} z^{-i+2} +
    A_{-i+3} z^{-i+3}  +\cdots \right) z^i\mbf{Y} \\
&\le z^i\left(B_{-i} z^{-i}\mbf Y+A^*(z)\mbf Y \right)
\le z^i\left(\alpha z^{-N}\mbf {e}^{t}+A^*(z)\mbf Y \right)\\
&= z^i\left(\alpha z^{-N}\mbf {e}^{t} +(1-\delta)\mbf Y\right)\le z^i\left(\alpha z^{-N} +(1-\delta)\right)\mbf Y\\
&\le (1-\frac\delta 2) \mbf v_i,
\endaligned
\end{equation}
where $  \mbf{v_i} = (V({i,1}), V({i,2}), \ldots, V({i,m}))^t=z^i\mbf Y$.

For any $0< i \le N$,
\[
    B_{-i} \mbf{v_0} + A_{-i+1} \mbf{v_1}  + A_{-i+2} \mbf{v_2} +
    A_{-i+3} \mbf{v_3}+
    \cdots \leq   \alpha\mbf {e}^{t} + z^N A^*(z)\mbf Y<\infty,
\]
and for $i=0$,
\[
    B_{0} \mbf{v_0}  + B_{1} \mbf{v_1} + B_{2} \mbf{v_2} +
    B_{3} \mbf{v_3}    +\cdots =B_0\mbf Y+B^*_+(z)\mbf Y<\infty.
\]
\S

The equivalence of geometric ergodicity and the geometric
stationary tail is interesting since results of ergodicity could
directly lead to results of the stationary tail asymptotics and
vice versa.

\section{Strong ergodicity}

In this section we show that the phase process $A$ is not stochastic if and only if
the $GI/G/1$-type Markov chain is strongly ergodic.

\begin{Thm}\label{thm:strong}
The $GI/G/1$-type Markov chain is strongly ergodic if and only if
 $A$ is not stochastic.
\end{Thm}

\Proof  Assume first that  $A$ is not stochastic.
Let $\tau_{0}=\inf\{n\ge0: X_n\in L_0\}$, then by Proposition 3.3 on page 216 in Anderson \cite{Andersonmc},
we only need to prove that there exists an $M<\infty$ such that
$$
x_{i,r}:=\mathbb E_{(i,r)}\tau_{0}\le M,\quad \text{for all}\  1\le i<\infty \ \text{and} \ 1\le r\le m.
$$
Let $\mbf X_i=(x_{i,1},\cdots,x_{i,m})^t$ for $i\ge1$ and $\mbf X_0=(x_{0,1},\cdots,x_{0,m_{0}})^t=\mbf 0$,
then $\{\mbf X_i: i\ge0\}$ is
the minimal
non-negative solution of
\begin{equation}\label{se1}
\mbf Y_i=\sum_{k\not=i}P_{i,k}\mbf Y_k+\mbf {e}^{t}, \quad i\ge1, \quad
\end{equation}
where $\mbf Y_0=\mbf 0$.

Since $A$ is irreducible and not stochastic,
$(I-A)^{-1}=I+A+A^2+\cdots$ exists and is finite. In the
following, we will prove that $ \mbf X_i\le (I-A)^{-1}\mbf e$. In
fact, by a standard procedure for the minimal non-negative
solution (for example, see \cite{hg}), set $\mbf X_i^{(0)}=0$, and
for $n\ge1$, let
\[
\mbf X_i^{(n)}=\sum_{k\not=i}P_{i,k}\mbf X_k^{(n-1)}+\mbf {e}^{t}, \quad i\ge1, \quad \text{and $\mbf X^{(n)}_0=0$,}
\]
we can then inductively prove that if $\mbf X_i^{(n-1)}\le
(I-A)^{-1}\mbf {e}^{t}$, then for $i\ge1$,
\[
\mbf X_i^{(n)} \le  \sum_{k\not=i}P_{i,k}(I-A)^{-1}\mbf {e}^{t}+\mbf {e}^{t} \le A(I-A)^{-1}\mbf {e}^{t}+\mbf {e}^{t}=(I-A)^{-1}\mbf {e}^{t}.
\]
Thus, $\mbf X_i=\lim_{n\rightarrow\infty}\mbf X_i^{(n)}\le (I-A)^{-1}\mbf {e}^{t}<\infty$.

For the converse, assume that $A$ is stochastic. It is obvious that
$P$ is a Feller transition matrix, thus $P$ cannot be strongly ergodic
by Proposition 2.3 in Hou and Liu \cite{HouEergodicity}.
\S


\section{Polynomial ergodicity}

 As indicated in the introduction, a sufficient condition for polynomial ergodicity was obtained in H{\o}jgaard and M{\o}ller~\cite{Hojgaard1996}. In this section,
 we prove that it is also necessary . In addition, we provide another necessary and sufficient condition for polynomial ergodicity. To achieve this goal, in the first sub-section, we construct a control function $h$ based on transition probabilities and provide a lower bound for the first hitting time. We also find a relationship between the first hitting times and
transition probabilities. In the second sub-section, proofs of the main results are provided.

\subsection{Lower bounds for the first hitting time}\label{sec:04}

In this sub-section, we derive lower bounds for the first hitting time $\tau_{i}=\inf\{ n\geq 1: (X_{n},Y_{n})\in L_{i} \}$, for any finite $i$, which are useful to discuss necessary conditions for ergodicity. We provide a lemma for $\tau_{0}$. Similar results can be obtained for any $i \neq 0$ by the same argument. In fact, we can have a corresponding result for the first hitting time of any finite set instead of a level. We also obtain a relationship between the first hitting time and transition probabilities.

\begin{Lem} \label{lemergodic01}
If $\sum\limits_{k=0}^{\infty} kB_{k}\be^{t} < \boldsymbol{\infty}^{t}$ and $\sum\limits_{k=0}^{\infty} kA_{k}\be^{t} < \boldsymbol{\infty}^{t}$, then (1) for $i$ large enough, we have,
$\p(C_{i}) > \frac{1}{2}$, where  $C_{i}=\left\{(X_{s},Y_{s})\in \left(\bigcup_{u=i}^{\lfloor i+\frac{i}{2}\rfloor} L_{u}\right)\quad \text{for} \quad s=0,1,\dots,\left\lfloor\frac{i}{4\mu}\right\rfloor\right\}$ and
$\lfloor x \rfloor$ denotes the largest integer equal to or smaller than $x$; (2) for each sample on the event to $C_{i}$ which satisfies (1),
\eqnb \label{seergodic0501}
 \tau_{0}-1 \geq \frac{i-1}{4\mu},
\eqne
where $\mu > 0$ is some fixed constant.
\end{Lem}
\Proof Since we consider the level independent $GI/G/1$-type Markov chain, we have, for $i\geq 1$ and $j=1,\dots,m$,
\[
       \p_{(i,j)}\left((X_{1},Y_{1})\in \left(\bigcup_{u=i}^{i+(k-1)} L_{u}\right)^{c}\right)
    =  \p_{(1,j)}\left((X_{1},Y_{1})\in \left(\bigcup_{u=1}^{1+(k-1)} L_{u}\right)^{c}\right), \quad \text{for any}\, k>0,
\]
where $L^{c}$ denotes the complement of $L$ and $\p_{(i,j)}((X_{1},Y_{1})\in \cdot\,\,)=\p((X_{1},Y_{1})\in \cdot\,\,|(X_{0},Y_{0})=(i,j))$.
Without loss of generality, we assume that
\begin{align*}
      & \p_{(0,j_{1})}\left((X_{1},Y_{1})\in \left(\bigcup_{u=0}^{k-1} L_{u}\right)^{c}\,\,\,\right) \\
    = & \max_{j}\left\{\p_{(0,j)}\left((X_{1},Y_{1})\in \left(\bigcup_{u=0}^{k-1} L_{u}\right)^{c}\,\,\,\right),\quad j=1,\dots,m\right\},
\end{align*}
and
\begin{align*}
      & \p_{(1,j_{2})}\left((X_{1},Y_{1})\in \left(\bigcup_{u=1}^{1+(k-1)} L_{u}\right)^{c}\,\,\,\right)\\
    = & \max_{j}\left\{\p_{(1,j)}\left((X_{1},Y_{1})\in \left(\bigcup_{u=1}^{1+(k-1)} L_{u}\right)^{c}\,\,\,\right),\quad j=1,\dots,m\right\}.
\end{align*}
Let
\[
 h(k)=\max\left\{\p_{(0,j_{1})}\left((X_{1},Y_{1})\in \left(\bigcup_{u=0}^{k-1} L_{u}\right)^{c}\,\,\right), \, \p_{(1,j_{2})}\left((X_{1},Y_{1})\in \left(\bigcup_{u=1}^{1+(k-1)} L_{u}\right)^{c}\,\,\right) \right\}.
\]
Notice that
\[
 \sum\limits_{k=1}^{\infty} \p_{(0,j_{1})}\left((X_{1},Y_{1})\in \left(\bigcup_{u=0}^{k-1} L_{u}\right)^{c}\,\,\,\right) <\infty,
\]
and
\[
 \sum\limits_{k=0}^{\infty} \p_{(1,j_{2})}\left((X_{1},Y_{1})\in \left(\bigcup_{u=1}^{1+(k-1)} L_{u}\right)^{c}\,\,\,\right) <\infty
\]
since $\sum\limits_{k=0}^{\infty} kB_{k}\be^{t} < \boldsymbol{\infty}^{t}$ and $\sum\limits_{k=0}^{\infty} kA_{k}\be^{t} < \boldsymbol{\infty}^{t}$. Hence
\[
 \sum\limits_{k=1}^{\infty} h(k) <\infty.
\]
It is easy to know that
for all $i\geq 0$, $ j=1,\dots,m$ and $k > 0$,
\[
  \p_{(i,j)}\left((X_{1},Y_{1})\in \left(\bigcup_{u=i}^{i+(k-1)} L_{u}\right)^{c}\,\,\,\right)\leq h(k),
\]
and $h(k)$ is a non-increasing function. Therefore there exists a sequence of $i.i.d.$ random variables $W_{n}>0$ with mean $\mu=\e[W_{n}]$ such that for all $i\geq 0$ and all $k\geq 0$,
\[
  \p_{(i,j)}\left((X_{1},Y_{1})\in \left(\bigcup_{u=i}^{i+(k-1)} L_{u}\right)^{c}\,\,\,\right)\leq \p(W_{n} \geq k).
\]
By the weak law of large numbers, we claim that for any $\epsilon>0$,
\[
 \lim_{n\to \infty}\p(S_{n}\leq (\mu+\epsilon)n)=1,
\]
where $S_{n}=W_{1}+\dots+W_{n}$. Hence there exists an $N$ large enough such that for $n\geq N$,
\[
 \p(S_{n}<2\mu n)\geq \frac{1}{2}.
\]
Then, a stochastic comparison argument yields, for all $i\geq 0$ and all $n\geq N$,
\eqnb \label{seergodic0504}
 \p_{(i,j)}\left((X_{s},Y_{s})\in \left(\bigcup_{u=i}^{\lfloor i+2\mu n\rfloor} L_{u}\right)\quad \text{for} \quad s=0,1,\dots n\right)\geq \frac{1}{2}.
\eqne
For $i$ large enough such that $\frac{i}{4\mu}\geq N$, we have from~(\ref{seergodic0504}) with $n=\left\lfloor\frac{i}{4\mu}\right\rfloor \geq N$ that
\eqnb \label{seergodic0505}
 \p_{(i,j)}\left((X_{s},Y_{s})\in \left(\bigcup_{u=i}^{\lfloor i+\frac{i}{2}\rfloor} L_{u}\right)\quad \text{for} \quad s=0,1,\dots,\left\lfloor\frac{i}{4\mu}\right\rfloor \right)\geq \frac{1}{2}. \notag
\eqne
Set $C_{i}=\left\{(X_{s},Y_{s})\in \left(\bigcup_{u=i}^{\lfloor i+\frac{i}{2}\rfloor} L_{u}\right)\quad \text{for} \quad s=0,1,\dots,\left\lfloor\frac{i}{4\mu}\right\rfloor\right\}$, then we have that  $\p(C_{i}) > \frac{1}{2}$.
Therefore, for $i$ large enough, we have, for each sample on the event $C_{i}$,
\[
 \tau_{0}-1 \geq \left\lfloor\frac{i}{4\mu}\right\rfloor \geq \frac{i-1}{4\mu}.
\]
\S

\begin{Cor} \label{corergodic03}
Given any non-negative, non-decreasing and measurable function $f$, if $\,\,\e_{(0,j)}[f(\tau_{0})]<\infty$ for\,\, $j=1,2\dots,m$,\,\, $\sum\limits_{k=0}^{\infty} kB_{k}\be^{t} < \boldsymbol{\infty}^{t}$, and  $\sum\limits_{k=0}^{\infty} kA_{k}\be^{t} < \boldsymbol{\infty}^{t}$, then for $i$ large enough,  we have
\eqnb \label{seergodic0516}
 \e_{(i,j)}[f(\tau_{0}-1)] \geq \frac{1}{2}f\left(\frac{i-1}{4\mu}\right),
\eqne
where $\mu > 0$ is some fixed constant.
\end{Cor}
\Proof From Lemma~\ref{lemergodic01}, we know that for $i$ large enough such that $\left\lfloor\frac{i}{4\mu}\right\rfloor\geq N$,
\eqnb \label{seergodic0517}
 \p_{(i,j)}\left((X_{s},Y_{s})\in \left(\bigcup_{u=i}^{\lfloor i+\frac{i}{2}\rfloor} L_{u}\right)\, \text{for} \,\,s=0,1,\dots,\left\lfloor\frac{i}{4\mu}\right\rfloor \right)\geq \frac{1}{2}.
\eqne Hence for $i$ large enough, we have for the above event
\[
 \tau_{0}-1 \geq \left\lfloor\frac{i}{4\mu}\right\rfloor \geq \frac{i-1}{4\mu},
\]
and therefore
\[
 f(\tau_{0}-1) \geq f\left(\left\lfloor\frac{i}{4\mu}\right\rfloor\right) \geq f\left(\frac{i-1}{4\mu}\right).
\]
For $i$ sufficiently large, this event has a probability of at
least $\frac{1}{2}$ by~(\ref{seergodic0517}) and therefore for $i$
sufficiently large, \eqnb
 \e_{(i,j)}[f(\tau_{0}-1)] \geq \frac{1}{2}f\left(\frac{i-1}{4\mu}\right), \notag
\eqne
where $\mu=\e[W_{n}]>0$ is some fixed constant.
\S

\begin{remark}
Specifically, if for $l\in\{2,3,\dots\}$, $f(x)=x^{l}$, then for
$i$ large enough,  we have \eqnb \label{seergodic0502}
 \e_{(i,j)}\left[(\tau_{0}-1)^{l}\right] \geq \frac{1}{2}\left(\frac{i-1}{4\mu}\right)^{l},
\eqne
where $\mu > 0$ is some fixed constant and  $\e_{(i,j)}\left[(\tau_{0}-1)^{l}\right]=\e\left[(\tau_{0}-1)^{l}|(X_{0},Y_{0})=(i,j)\right]$.
\end{remark}

Next, we discuss the relationship between the first hitting time
and one-step transition probabilities, which leads to a necessary
condition for polynomial ergodicity.
\begin{Lem} \label{lemergodic02}
If for $l\in\{2,3,\dots\}$,  $\e_{(0,j)}[\tau_{0}^{l}] < \infty$, for $j=1,\dots,m$, we have
\eqnb \label{seergodic0503}
 \e_{(0,j)}\left[\tau_{0}^{l}\right] = \sum\limits_{k=0}^{\infty}\sum\limits_{i=1}^{m} \e_{(k,i)}\left[(\tau_{0}+1)^{l}\right] \p_{(0,j)}((X_{1},Y_{1})=(k,i)),
\eqne
and
\eqnb \label{seergodic0518}
 \e_{(1,j)}\left[\tau_{0}^{l}\right] = \sum\limits_{k=0}^{\infty}\sum\limits_{i=1}^{m} \e_{(k,i)}\left[(\tau_{0}+1)^{l}\right] \p_{(1,j)}((X_{1},Y_{1})=(k,i)),
\eqne
where $\p_{(0,j)}((X_{1},Y_{1})=(k,i))=\p((X_{1},Y_{1})=(k,i)|(X_{0},Y_{0})=(0,j))$ and\, $\p_{(1,j)}((X_{1},Y_{1})=(k,i))=\p((X_{1},Y_{1})=(k,i)|(X_{0},Y_{0})=(1,j))$.
\end{Lem}
\Proof Given $l\in\{2,3,\dots\}$, we have that for  $j=1,\dots,m$,

\begin{align} \label{seergodic0508}
  & \e_{(0,j)}\left[\tau_{0}^{l}\right] \notag\\
= & \e_{(0,j)}\left[\e\left[\tau_{0}^{l}|(X_{1},Y_{1})\right]\right] \notag\\
= & \sum\limits_{k=0}^{\infty}\sum\limits_{i=1}^{m} \e\left[\tau_{0}^{l}|(X_{1},Y_{1})=(k,i)\right] \p_{(0,j)}((X_{1},Y_{1})=(k,i)) \notag\\
= & \sum\limits_{k=0}^{\infty}\sum\limits_{i=1}^{m}\sum\limits_{n=1}^{\infty} n^{l}\p(\tau_{0}=n|(X_{1},Y_{1})=(k,i))
\p_{(0,j)}((X_{1},Y_{1})=(k,i)) \notag\\
= & \sum\limits_{k=0}^{\infty}\sum\limits_{i=1}^{m}\sum\limits_{n=1}^{\infty} n^{l}\p(\tau_{0}=n-1|(X_{0},Y_{0})=(k,i)) \p_{(0,j)}((X_{1},Y_{1})=(k,i)) \notag\\
= & \sum\limits_{k=0}^{\infty}\sum\limits_{i=1}^{m}\sum\limits_{n=1}^{\infty} (n+1)^{l}\p(\tau_{0}=n|(X_{0},Y_{0})=(k,i)) \p_{(0,j)}((X_{1},Y_{1})=(k,i)) \notag\\
= & \sum\limits_{k=0}^{\infty}\sum\limits_{i=1}^{m} \e_{(k,i)}\left[(\tau_{0}+1)^{l}\right] \p_{(0,j)}((X_{1},Y_{1})=(k,i)). \notag
\end{align}
Similarly, we can prove the other case.
\S

\begin{remark}
In fact, for any non-negative and measurable function $f$, we can have decompositions corresponding to Lemma~\ref{lemergodic02} for $f(\tau_{0})$.
\end{remark}

\subsection{Necessary and sufficient conditions for polynomial ergodicity}\label{sec:002}

In this sub-section,  we discuss polynomial ergodicity for the $GI/G/1$-type Markov chain and provide two necessary and sufficient conditions.

The $GI/G/1$-type Markov chain is called polynomial ergodic of
degree $l$ for $l\in\{1, 2, \dots\}$, if for all $j=1,\dots,m$,
\eqnb\label{poly1}
\e_{(0,j)}\left[\tau_{\boldsymbol{0}}^{l}\right]<\infty. \eqne If
$l=1$ in equation~\ref{poly1}, this definition coincides with the
ordinary ergodicity.

We first use the lemmas from the previous sub-section to obtain a
necessary condition for polynomial ergodicity.

\begin{Thm} \label{thmergodic01}
If for $l\in\{2,3,\dots\}$, $\e_{(0,j)}\left[\tau_{0}^{l}\right]<\infty $ for all $j=1,\dots,m$, and $\sum\limits_{k=0}^{\infty} kA_{k}\be^{t} <\boldsymbol{\infty}^{t}$, then
\[
 \sum_{k} k^{l}A_{k}\be^{t} < \boldsymbol{\infty}^{t},
\]
and
\[
 \sum_{k} k^{l}B_{k}\be^{t} < \boldsymbol{\infty}^{t}.
\]
\end{Thm}

\Proof Given $l\in\{2,3,\dots\}$, by Corollary~\ref{corergodic03}, for $j=1,\dots,m$, we have for all $k>N$, where $N$ is large enough,
\eqnb \label{seergodic0509}
 \e_{(k,j)}\left[(\tau_{0}-1)^{l}\right] \geq \frac{1}{2}\left(\frac{k-1}{4\mu}\right)^{l}.
\eqne
Then, by Lemma~\ref{lemergodic02},
\begin{align} \label{seergodic0510}
   \e_{(0,j)}\left[\tau_{0}^{l}\right]
= & \sum\limits_{k=0}^{\infty}\sum\limits_{i=1}^{m} \e_{(k,i)}\left[(\tau_{0}+1)^{l}\right] \p_{(0,j)}((X_{1},Y_{1})=(k,i)) \notag\\
= & \sum\limits_{k=0}^{N}\sum\limits_{i=1}^{m} \e_{(k,i)}\left[(\tau_{0}+1)^{l}\right] \p_{(0,j)}((X_{1},Y_{1})=(k,i)) \notag\\
 & \,\,\,\,\,\,\,\,\,\,\,\, + \sum\limits_{k=N+1}^{\infty}\sum\limits_{i=1}^{m} \e_{(k,i)}\left[(\tau_{0}+1)^{l}\right] \p_{(0,j)}((X_{1},Y_{1})=(k,i)). \notag\\
\geq & \left[\sum\limits_{k=N+1}^{\infty} k^{l}B_{k}\be^{t}\right]_{j},\notag
\end{align}
where $\left[\sum\limits_{k=N+1}^{\infty} k^{l}B_{k}\be^{t}\right]_{j}$ denotes the $j$th element of vector $\sum\limits_{k=N+1}^{\infty} k^{l}B_{k}\be^{t}$.
Hence
\[
 \sum\limits_{k=N+1}^{\infty} k^{l}B_{k}\be^{t} < \boldsymbol{\infty}^{t},
\]
and therefore
\[
 \sum_{k} k^{l}B_{k}\be^{t} < \boldsymbol{\infty}^{t}.
\]
Similarly, we can obtain
\[
 \sum_{k} k^{l}A_{k}\be^{t} < \boldsymbol{\infty}^{t}.
\]
\S

\begin{remark}
This theorem can be extended to a class of more general non-negative and non-decreasing rate functions.
\end{remark}

A necessary and sufficient condition for polynomial ergodicity can
now be obtained since H{\o}jgaard and
M{\o}ller~\cite{Hojgaard1996} have already showed that the
converse of Theorem~\ref{thmergodic01} is also true.  We state it
as follows.
\begin{Thm} \label{thmergodic02}
If $\sum\limits_{k=0}^{\infty} kA_{k}\be^{t} <\boldsymbol{\infty}^{t}$, then for $l\in\{2,3,\dots\}$, $\e_{(0,j)}\left[\tau_{0}^{l}\right]<\infty$ for all $j=1,\dots,m$ if and only if $\sum_{k} k^{l}A_{k}\be^{t} < \boldsymbol{\infty}^{t}$ and $\sum_{k} k^{l}B_{k}\be^{t} < \boldsymbol{\infty}^{t}$.
\end{Thm}
\Proof For necessity, it follows from Theorem~\ref{thmergodi02}.
\S

Next, we provide another necessary and sufficient condition for polynomial ergodicity based on the relationship between ergodicity and the tail behavior of the stationary distribution for the $GI/G/1$-type Markov chain. The following lemma from Jarner and Tweedie~\cite{Jarner2001} is needed.
\begin{Lem} \label{lemergodic03}
Assume that for $l\in\{2,3,\dots\}$, $\e_{(0,j)}\left[\tau_{0}^{l}\right]<\infty\,$ for $j=1,\dots,m$, there exists a finite set (without loss of generality, we assume that this finite set is $L_{0}$) such that
\eqnb \label{seergodic0511}
 \sum\limits_{i=0}^{\infty} \e_{(i,j)}\left[\tau_{0}^{l}\right]\pi_{ij} < \infty,
\eqne
where $\pi_{ij}$ is the $j$th element of $\bpi_{i}$.
\end{Lem}

Using Lemma~\ref{lemergodic02}, Lemma~\ref{lemergodic03} and
factorization results, we have the following condition.
\begin{Thm} \label{thmergodi02}
If $\sum\limits_{k=0}^{\infty} kA_{k}\be^{t} <\boldsymbol{\infty}$, then for $l\in\{2,3,\dots\}$, $\e_{(0,j)}\left[\tau_{0}^{l}\right]<\infty\,$ for all $j=1,\dots,m$ if and only if $\sum_{i} i^{l}\pi_{ij} <\infty $.
\end{Thm}

\Proof For necessity, by Corollary~\ref{corergodic03}, for a given $l\in\{2,3,\dots\}$ and all $j=1,\dots,m$, we have for all $i>N$, where $N$ is large enough,
\eqnb
 \e_{(i,j)}\left[(\tau_{0}-1)^{l}\right] \geq \frac{1}{2}\left(\frac{i-1}{4\mu}\right)^{l}.
\eqne
Then by Lemma~\ref{lemergodic03},
\begin{align} \label{seergodic0510}
 \infty> & \sum\limits_{i=0}^{\infty} \e_{(i,j)}\left[\tau_{0}^{l}\right]\pi_{ij} \notag\\
= & \sum\limits_{i=0}^{N} \e_{(i,j)}\left[\tau_{0}^{l}\right]\pi_{ij} + \sum\limits_{i=N+1}^{\infty} \e_{(i,j)}\left[\tau_{0}^{l}\right]\pi_{ij} \notag\\
\geq & \sum\limits_{i=N+1}^{\infty} i^{l}\pi_{ij}.\notag
\end{align}
Thus
\[
 \sum_{i} i^{l}\pi_{ij} <\infty.
\]

For sufficiency, we use factorization results and generating function techniques in our analysis.
From~(\ref{SE:0103}) we have that for $ 0< z <1$,
\eqnb\label{seergodic0512}
   \bpi^{*}(z)[I - R^{*}(z)] = \bpi_{0}R^{*}_{0}(z),\notag
\eqne
and thus
\eqnb\label{seergodic0513}
   \bpi^{*}(z) = \bpi^{*}_{0}R^{*}_{0}(z)[I - R^{*}(z)]^{-1},
\eqne
since $I - R(z)$ is invertible. Taking the $l$th $(l=2,3,\dots)$ derivative on the both sides of equation~(\ref{seergodic0513}), we have
\begin{align} \label{seergodic0514}
   &\sum\limits_{k=l-1}^{\infty} k(k-1)\cdots(k-l+1)z^{k-l}\bpi_{k}  \notag\\
 = & \bpi_{0}\sum\limits_{k=l-1}^{\infty} k(k-1)\cdots(k-l+1)z^{k-l} R_{0, k} \sum\limits_{k=1}^{\infty} z^{k}\sum\limits_{n=0}^{\infty} R_{k}^{n*}  \\
   &\hspace{15mm} + \bpi_{0}\sum\limits_{k=1}^{\infty} z^{k}R_{0, k}\sum\limits_{k=l-1}^{\infty} k(k-1)\cdots(k-l+1)z^{k-l}\sum\limits_{n=0}^{\infty} R_{k}^{n*}
   + \boldsymbol{c}_{l-1}(z) \notag ,
\end{align}
where $\boldsymbol{c}_{l-1}(z)$ is the summation of all terms of the form
\[
 \bpi_{0}\sum\limits_{k=p-1}^{\infty} k(k-1)\cdots(k-p+1)z^{k-p}R_{0, k} \sum\limits_{k=q-1}^{\infty} k(k-1)\cdots(k-q+1)z^{k-q}\sum\limits_{n=0}^{\infty} R_{k}^{n*}
\]
and
\[
 p+q \leq l-1.
\]
There are only finitely many such terms of this kind. Let $z\to 1-$ in~(\ref{seergodic0514}), we obtain
\begin{align} \label{seergodic0515}
   &\sum\limits_{k=l-1}^{\infty} k(k-1)\cdots(k-l+1)\bpi_{k}  \notag\\
 = & \bpi_{0}\sum\limits_{k=l-1}^{\infty} k(k-1)\cdots(k-l+1) R_{0, k} \sum\limits_{k=1}^{\infty} \sum\limits_{n=0}^{\infty} R_{k}^{n*}  \notag\\
   &\hspace{15mm} + \bpi_{0}\sum\limits_{k=1}^{\infty} R_{0, k}\sum\limits_{k=l-1}^{\infty} k(k-1)\cdots(k-l+1)\sum\limits_{n=0}^{\infty} R_{k}^{n*}
    + \boldsymbol{c}_{l-1}(1) \notag .
\end{align}
Since $\sum_{i} i^{l}\bpi_{i} <\boldsymbol{\infty} $ and
\[
\lim_{k\to \infty}\frac{k(k-1)\cdots(k-l+1)}{k^{l}}=1,
\]
it follows that
\[
 \sum\limits_{k=l-1}^{\infty} k(k-1)\cdots(k-l+1)\bpi_{k} < \infty.
\]
Therefore
\[
 \sum\limits_{k=l-1}^{\infty} k(k-1)\cdots(k-l+1) R_{0, k} < \infty
\]
and
\[
 \sum\limits_{k=l-1}^{\infty} k(k-1)\cdots(k-l+1)\sum\limits_{n=0}^{\infty} R_{k}^{n*} < \infty.
\]
Thus
\[
 \sum\limits_{k=1}^{\infty} k^{l} R_{0, k} < \infty
\]
and
\[
 \sum\limits_{k=1}^{\infty} k^{l}\sum\limits_{n=0}^{\infty} R_{k}^{n*} < \infty.
\]
From Theorem~1 in~\cite{Grassmann:14} and Theorem~12 in~\cite{Zhao:26}, we have for $k\geq 1$,
$ R_{0,k} \geq B_{k} $ and $ R_{k} \geq A_{k} $. Hence,
\[
 \sum\limits_{k=1}^{\infty} k^{l} B_{k}\be^{t} \leq \sum\limits_{k=1}^{\infty} k^{l} R_{0, k}\be^{t} < \boldsymbol{\infty}^{t}
\]
and
\[
 \sum\limits_{k=1}^{\infty} k^{l} A_{k}\be^{t} \leq \sum\limits_{k=1}^{\infty} k^{l}\sum\limits_{n=0}^{\infty} A_{k}^{n*}\be^{t} \leq \sum\limits_{k=1}^{\infty} k^{l}\sum\limits_{n=0}^{\infty} R_{k}^{n*}\be^{t} < \boldsymbol{\infty}^{t}.
\]
Finally, by Theorem~\ref{thmergodic02}, we know that $\e_{(0,j)}\left[\tau_{0}^{l}\right]<\infty\,$.
\S

\end{document}